\documentclass[letterpaper, 10 pt, conference]{ieeeconf}  

\IEEEoverridecommandlockouts                              

\usepackage{amsmath,amssymb,amsfonts}
\usepackage{algorithmic}
\usepackage{graphicx,graphics,wrapfig,lipsum}
\usepackage{textcomp}
\usepackage{xcolor}

\makeatletter
\def\maketag@@@#1{\hbox{\m@th\normalfont\normalsize#1}}
\makeatother

\usepackage{mathtools}
\usepackage{bm}
\usepackage{array}
\usepackage{mathrsfs}

\newlength\fwidth

\usepackage{ifthen}
\usepackage{xargs}
\usepackage{pdfpages}
\usepackage{calc}
\usepackage [english]{babel}

\usepackage{microtype}
\usepackage{svg}

\allowdisplaybreaks

\newcolumntype{P}[1]{>{\centering\arraybackslash}p{#1}}
\newcolumntype{M}[1]{>{\centering\arraybackslash}m{#1}}
\newcolumntype{B}[1]{>{\centering\arraybackslash}b{#1}}

\newcommand\myHLine[1]{\noalign{\hrule height #1}}

\title{\LARGE \bf
	Control for autonomous vehicles in collision avoidance maneuvers : LPV modeling and static feedback controller
}

\author{Dario Penco$^{1,}$$^{2}$, Joan Davins-Valldaura$^{2}$, Emmanuel Godoy$^{1}$, Pedro Kvieska$^{2}$ and Giorgio Valmorbida$^{1}$
\thanks{$^{1}$Université Paris-Saclay, CNRS, CentraleSupélec,
Laboratoire des signaux et systèmes, 91190 Gif-sur-Yvette, France.
{\tt\small \{dario.penco, emmanuel.godoy, giorgio.valmorbida\}@centralesupelec.fr}}%
\thanks{$^{2}$Groupe Renault, Technocentre, 78288 Guyancourt, France.
{\tt\small \{dario.penco, joan.davins-valldaura, pedro.kvieska\}@renault.com}}%
}

\begin{document}
	
	\maketitle
	\thispagestyle{empty}
	\pagestyle{empty}
	
	\begin{abstract}
		
		This article presents a state feedback control design strategy for the stabilization of a vehicle along a reference collision avoidance maneuver.
		The stabilization of the vehicle is achieved through a combination of steering, acceleration and braking.
		A Linear Parameter-Varying (LPV) model is obtained from the linearization of a non-linear model along the reference trajectory.
		A robust state feedback control law is computed for the LPV model.
		Finally, simulation results illustrate the stabilization of the vehicle along the reference trajectory.

	\end{abstract}

	\section{Introduction}
	One of the current challenges in the development of autonomous driving is to guarantee safety for increased levels of autonomous driving \cite{sae_international_taxonomy_2018}.
	For low velocity maneuvers, such as parking, the trajectory planning and execution can rely on simple mathematical models.
	On the other hand, collision avoidance maneuvers, carried out in high velocities require more sophisticated control strategies to be automated and, due to its dynamics involving suspension, brakes and tire forces call for more complex models.\\
	In the context of autonomous driving, maneuvers with a maximum vehicle's acceleration of $0.3$g on just one vehicle axis are considered as low dynamics maneuvers.
	Maneuvers with an acceleration larger than $0.3$g and possibly a combination of longitudinal and lateral acceleration are considered \emph{high dynamics} maneuvers.\\
	Literature is rich of solutions for vehicle control in low dynamics maneuvers.
	In low dynamics the longitudinal and lateral dynamics of the vehicle are commonly considered uncoupled.
	The control of the longitudinal and lateral dynamics is hence treated separately.
	A common solution employed for longitudinal control is the \textit{Linear Proportional Derivative} (PID) control law \cite{hima_trajectory_2011}, \cite{ioannou_intelligent_1993}, \cite{naus_string-stable_2010} and \cite{persson_stop_1999}.
	Lateral control is usually based on the \textit{linear bicycle model} \cite{rajamani_vehicle_2012}.
	Several solutions use $H_2$ and $H_{\infty}$ control laws \cite{hima_trajectory_2011}, \cite{obrien_vehicle_1996} and \cite{shimakage_design_2002}. \\
	For high dynamics maneuvers the coupling of the vehicle's longitudinal and lateral dynamics cannot be neglected.
	Moreover in these conditions the combined longitudinal and lateral tire slip is an important factor related to the tire forces saturation \cite{michelin_pneu_2001}. \\
	Solutions for high dynamics strategies may also consider a decoupling of longitudinal and lateral dynamics \cite{kritayakirana_autonomous_2012}.
	However, since the cross effects on lateral and longitudinal variables are more significant, these strategies demand a finer tuning of the controller and a better knowledge of the vehicle parameters.
	These strategies may thus lack robustness. \\
	The use of joint lateral and longitudinal models will therefore be sought for high dynamics. \\
	\textit{Model Predictive Control} (MPC)  is a control strategy that attacks simultaneously the path generation and its control.
	Moreover, it allows to account for input- and state-constraints, which is useful since the obstacles and lane limits can be easily translated into constraints of the underlying optimization problem.
	Because of that in literature we find several solutions for vehicle control in high dynamics maneuvers based on the MPC control strategy.
	In \cite{falcone_predictive_2007} an MPC controller is developed based on a non-linear bicycle model that takes into consideration both longitudinal and lateral dynamics.
	The simulation and experimental tests show good results, but the computational load required by the controller is high.
	To solve this problem another MPC controller based on local linearizations of the non-linear model is proposed.
	The computational load is reduced, but the commands computed show chattering.
	In \cite{falcone_mpc-based_2008} an improvement is proposed.
	The controller is based on a four-wheel non-linear model.
	Also in this case the computational is considerable.
	A similar solution is proposed in \cite{yin_integration_2015}.
	In \cite{gao_tube-based_2014} \textit{tube-based} MPC is used to address also robustness criteria.
	In this case, to reduce the computational load by reducing the complexity of the model on which the controller is based, the longitudinal and lateral tire forces have been considered as the system's inputs.
	A low level controller is then required to compute the steering angle and wheel torques.
	Despite the good results present in literature concerning MPC, for commercial automobiles, the certification of a computationally demanding strategy may be more difficult.
	Due to the limited computational power the online optimization-based strategies as MPC may require ad hoc circuits. \\
	The goal of this paper is to propose a controller for vehicle stabilization along a reference high dynamics maneuver.
	To achieve this goal we propose a \textit{Linear Parameter-Varying} (LPV) model for the vehicle dynamics along the maneuver.
	The LPV framework offers methods for the synthesis of controllers that do not require a large computation load \cite{apkarian_self-scheduled_1995}, \cite{wu_control_1995}, \cite{biannic_commande_1996}, \cite{blanchini_set-theoretic_2015}. \\
	We develop the LPV model based on a nonlinear model that considers the coupling between longitudinal and lateral dynamics.
	Moreover, since in a high dynamics maneuvers the tires may operate on their physical limits, presenting combined tire slip, we further complexify the model by adding the tire nonlinear behavior.
	Clearly, other dynamics can also influence high maneuvers, such as the suspension and braking systems.
	In this paper we will not consider these dynamics as we focus on the study of the influence of the tire forces.\\
	The linearization of the proposed nonlinear model along a reference maneuver trajectory yields the LPV system we use for the synthesis of the control law.
	In particular, due to an algebraic loop introduced by the saturation function of the tires, we propose a strategy to obtain an uncertain parameter related to the slope of sector non-linearities. \\
	
	The outline of the article is:
	\begin{itemize}
		\item in Section \ref{sec_devMe_preliminaries} the non-linear vehicle model and the procedure to obtain the LPV model are discussed;
		\item the method used for the synthesis of the controller is shown in Section \ref{sec_devMe_main};
		\item in Section \ref{sec_devMe_numericalResults} the results of the application of the controller in a simulation environment are discussed;
		\item the conclusions and future work are discussed in Section \ref{sec_devMe_conclusion}.
	\end{itemize}

	\section{Preliminary results} \label{sec_devMe_preliminaries}
	
	\subsection{Vehicle model description and equations of motion} \label{subsec_devMe_nonLinearModel}
\begin{figure} [h]
		\centering
		\includegraphics[width = 1\columnwidth]{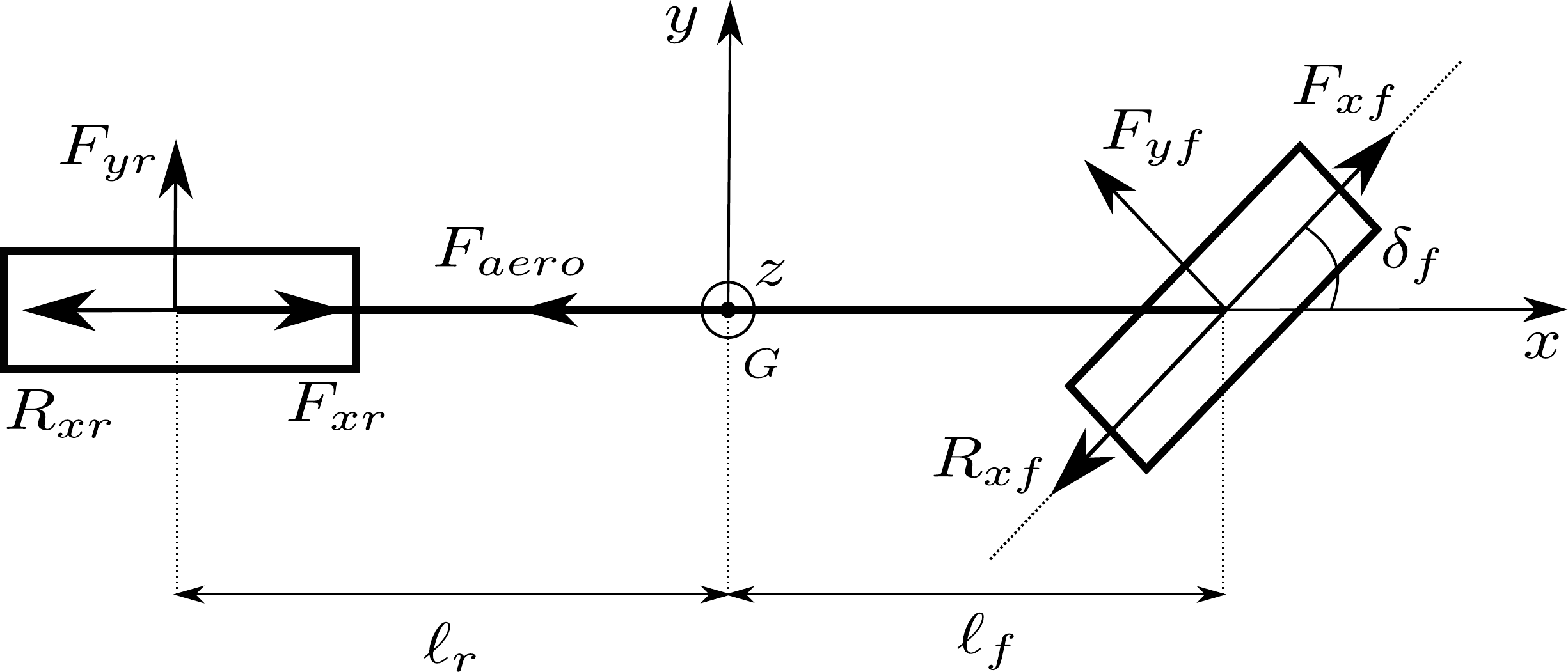}
		\caption{Vehicle's $xy$ plane, with the forces listed in Table \ref{table_devMe_forces}.}\label{fig_devMe_bicycle_from_above_noPitchRoll}
	\end{figure}

	To obtain a nonlinear model of the vehicle we consider that the two wheels of each axle have been regrouped at the center of the axle track.
	The vehicle degrees of freedom allow it to move on the $xy$ plane, so it can translate along the $x$ and $y$ axis and rotate around the $z$ axis.
	The vehicle axis are shown in Figure \ref{fig_devMe_bicycle_from_above_noPitchRoll}.\\
	In Table \ref{table_devMe_modelStateInputs} there are the model state variables and inputs. \\
	\begin{figure} [h]
		\centering
		\includegraphics[width = 0.55\columnwidth]{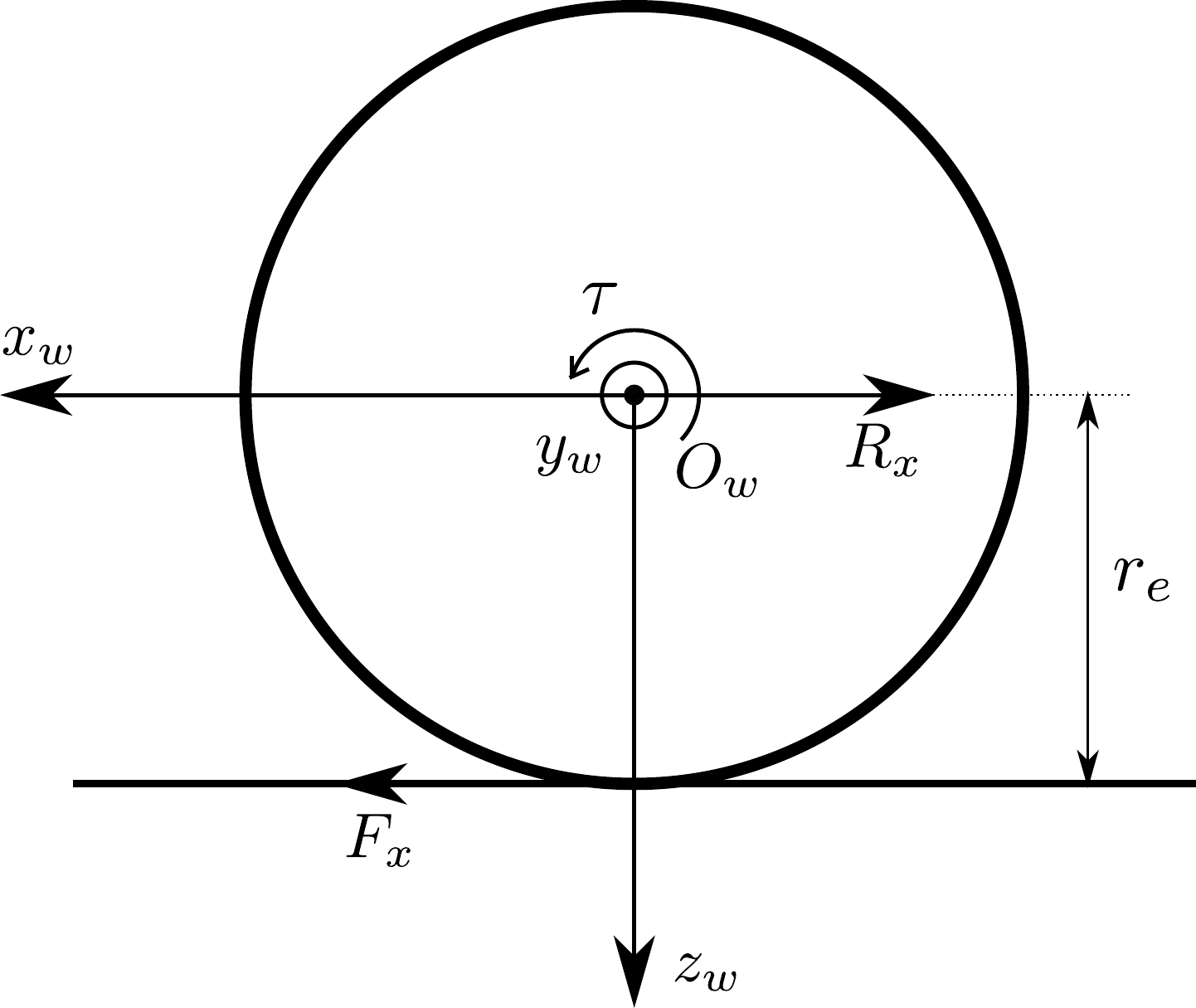}
		\caption{Wheel's $x_w z_w$ plane.}\label{fig_devMe_wheel_from_left_forces}
	\end{figure}
	\begin{table}[h!]
		\small
		\centering \renewcommand{\arraystretch}{1.5}
		\caption{Model state variables and inputs}
		\begin{tabular}{!{\vrule width 1pt} c c c !{\vrule width 1pt}}
			\myHLine{1pt}
			\multicolumn{3}{!{\vrule width 1pt} c !{\vrule width 1pt}}{State variables} \\
			\myHLine{1pt}
			$v$ & $\left[ \frac{\text{m}}{\text{s}}\right]$ & longitudinal speed\\
			\myHLine{0.1pt}
			$u$ & $\left[ \frac{\text{m}}{\text{s}}\right]$ & lateral speed\\
			\myHLine{0.1pt}
			$r$ & $\left[ \frac{\text{rad}}{\text{s}}\right]$ & yaw rate\\
			\myHLine{0.1pt}
			$\omega_{wf}$ & $\left[ \frac{\text{rad}}{\text{s}}\right]$ & front wheel rotational speed\\
			\myHLine{0.1pt}
			$\omega_{wr}$ & $\left[ \frac{\text{rad}}{\text{s}}\right]$ & rear wheel rotational speed\\
			\myHLine{0.1pt}
			$x$ & $\left[ \text{m} \right]$ & absolute longitudinal position\\
			\myHLine{0.1pt}
			$y$ & $\left[ \text{m} \right]$ & absolute lateral position\\
			\myHLine{0.1pt}
			$\psi$ & $\left[ \text{deg} \right]$ & yaw angle\\
			\myHLine{1pt}
			\noalign{\vskip 2pt}
			\myHLine{1pt}
			\multicolumn{3}{!{\vrule width 1pt} c!{\vrule width 1pt}}{System's inputs}\\
			\myHLine{1pt}
			$\delta_{f}$ & $\left[ \text{deg}\right]$ & front wheel steering angle \\
			\myHLine{0.1pt}
			$\tau_{wf}$ & $\left[ \text{Nm}\right]$ & front wheel torque \\
			\myHLine{0.1pt}
			$\tau_{wr}$ & $\left[ \text{Nm}\right]$ & rear wheel torque\\
			\myHLine{1pt}
		\end{tabular}
		\label{table_devMe_modelStateInputs}
	\end{table}
	\begin{table}[h!]
		\scriptsize
		\centering \renewcommand{\arraystretch}{1.5}
		\caption{Model's forces}
		\begin{tabular}{!{\vrule width 1pt} c c !{\vrule width 1pt} c c !{\vrule width 1pt}}
			\myHLine{1pt}
			$F_{xf}$ & front tire longitudinal force & $F_{xr}$ & rear tire longitudinal force \\
			\myHLine{0.1pt}
			$F_{yf}$ & front tire lateral force & $F_{yr}$ & rear tire lateral force \\
			\myHLine{0.1pt}
			$R_{xf}$ & front wheel rolling friction & $R_{xr}$ & rear wheel rolling friction \\
			\myHLine{0.1pt}
			$F_{aero}$ & aerodynamic drag & $P$ & weight \\
			\myHLine{0.1pt}
			$N_{f}$ & front wheel normal force & $N_{r}$ & rear wheel normal force \\
			\myHLine{0.1pt}
			\myHLine{1pt}
		\end{tabular}
		\label{table_devMe_forces}
	\end{table}
	~\\
	Table \ref{table_devMe_forces} contains the forces shown in Figures \ref{fig_devMe_bicycle_from_above_noPitchRoll} and \ref{fig_devMe_wheel_from_left_forces}.
	The forces applied on the front wheel have each a component on the $x$ and $y$ axis of the vehicle, that depends on $\delta_{f}$.
	We denote the components of $F_{xf}$ on the $xy$ plane:
	\begin{equation*}
	\begin{aligned}
		F_{xxf} &= F_{xf} \, \cos \delta_{f}\\
		F_{yxf} &= F_{xf} \, \sin \delta_{f}.
	\end{aligned}
	\end{equation*}
	The model equations of motion in the vehicle frame are:
	\small\begin{subequations} \label{eq_devMe_bicycle_EquationMotion_3DOF_vehicle}
	\begin{align}
	\begin{split}
	m \left(\dot{v} - r \, u\right)  &= 2 \, \left( F_{xxf} + F_{xr} - R_{xxf} - R_{xr} \right. \\
	& \left. - F_{xyf} \right) - F_{aero}
	\end{split} \\[8pt]
	m \left(\dot{u} + r \, v\right) &= 2 \, \left( F_{yxf} - R_{yxf} + F_{yyf} + F_{yr} \right) \\[8pt]
	I_{zz} \dot{r} &= 2 \, \left(F_{yxf} + F_{yyf} - R_{yxf}\right) \, \ell_f - 2 \,F_{yr} \, \ell_r \\[8pt]
	I_{w \, y} \, \dot{\omega}_{wf} &= \tau_{wf} - 2 \, F_{xf}\, r_{e}\\[8pt]
	I_{w \, y} \, \dot{\omega}_{wr} &= \tau_{wr} - 2 \, F_{xr}\, r_{e}.
	\end{align}
	\end{subequations}\normalsize
	For the sake of brevity, just the tires forces and the normal forces will be detailed.
	For the other forces description see \cite{rajamani_vehicle_2012} and \cite{jazar_vehicle_2017}.\\
	The vehicle position and orientation on the inertial frame are governed by:
	\begin{subequations} \label{eq_devMe_bicycle_posDynamics}
		\begin{align}
		\dot{x}  &= v \, \cos \psi - u \, \sin \psi \label{eq_devMe_bicycle_posDynamics_1}\\
		\dot{y}  &= u \, \cos \psi + v \, \sin \psi \label{eq_devMe_bicycle_posDynamics_2}\\
		\dot{\psi} &= r.
		\end{align}
	\end{subequations}
	
	\subsection{Tires forces} \label{subsec_devMe_tireModel}
\begin{figure} [h]
		\centering
		\includegraphics[width = 0.55\columnwidth]{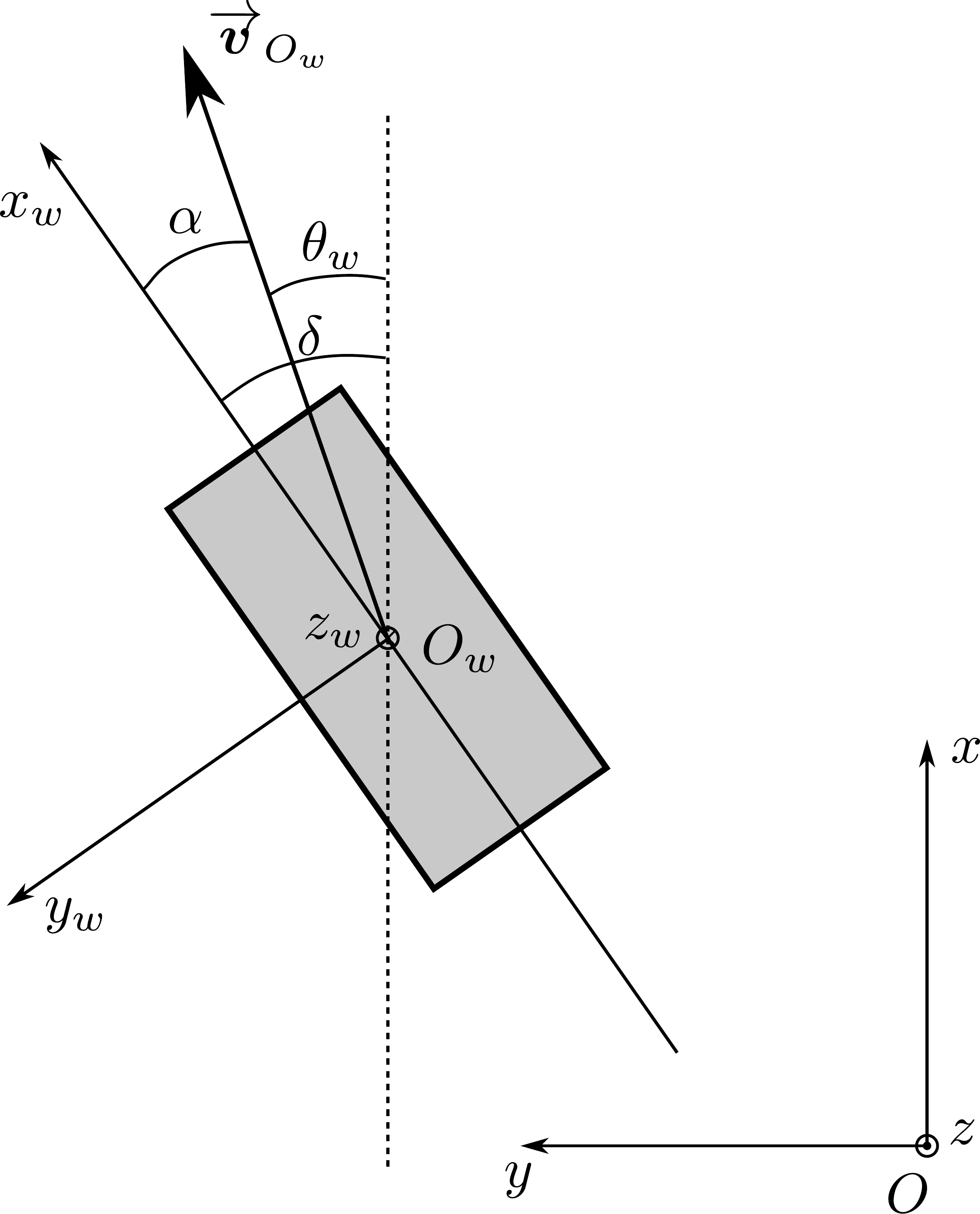}
		\caption{Wheel lateral slip angle $\alpha$.}\label{fig_devMe_wheel_slipAngle}
	\end{figure}

\begin{figure}[b]
		\centering
		\includegraphics[width=0.85\columnwidth]{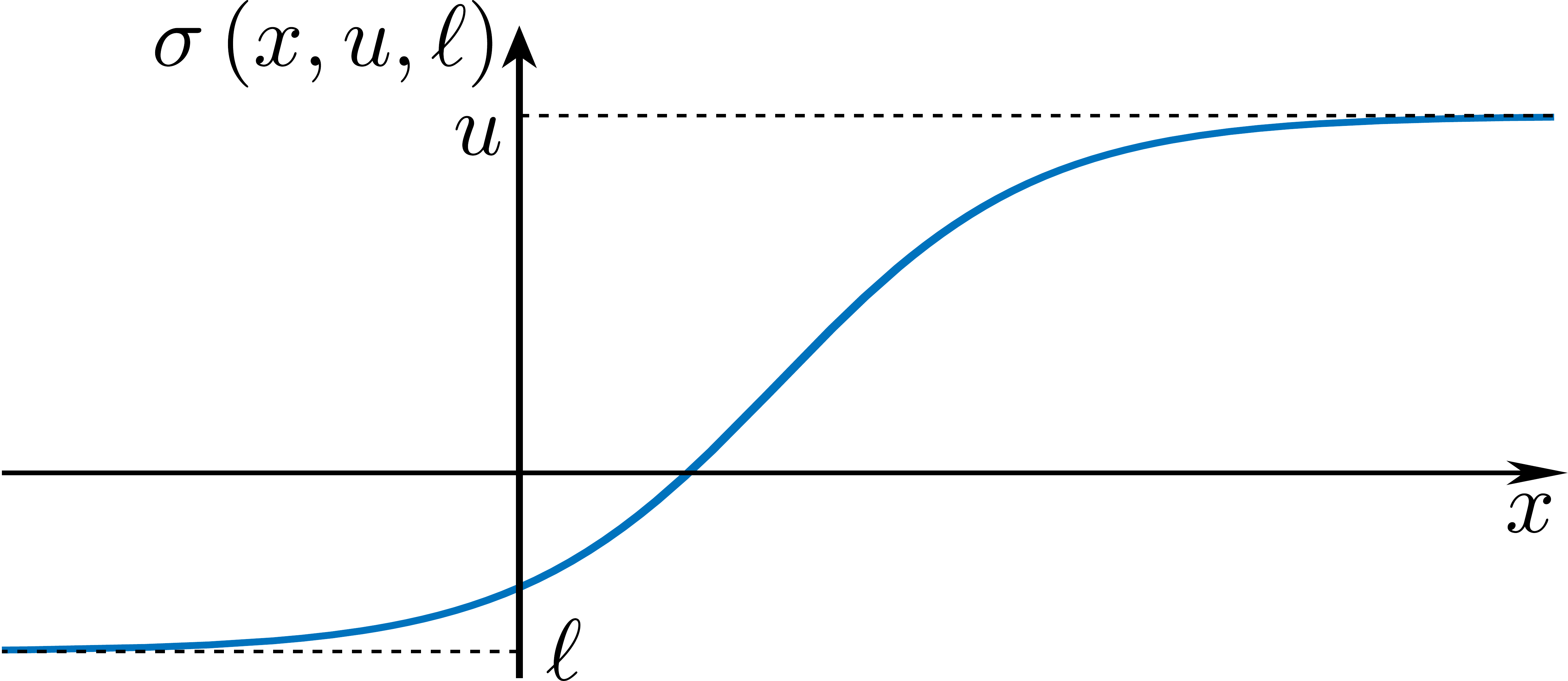}
		\caption{Logistic function.}\label{fig_devMe_logistic_function}
	\end{figure}

\begin{figure}[h!]
		\centering
		\includegraphics[width=0.85\columnwidth]{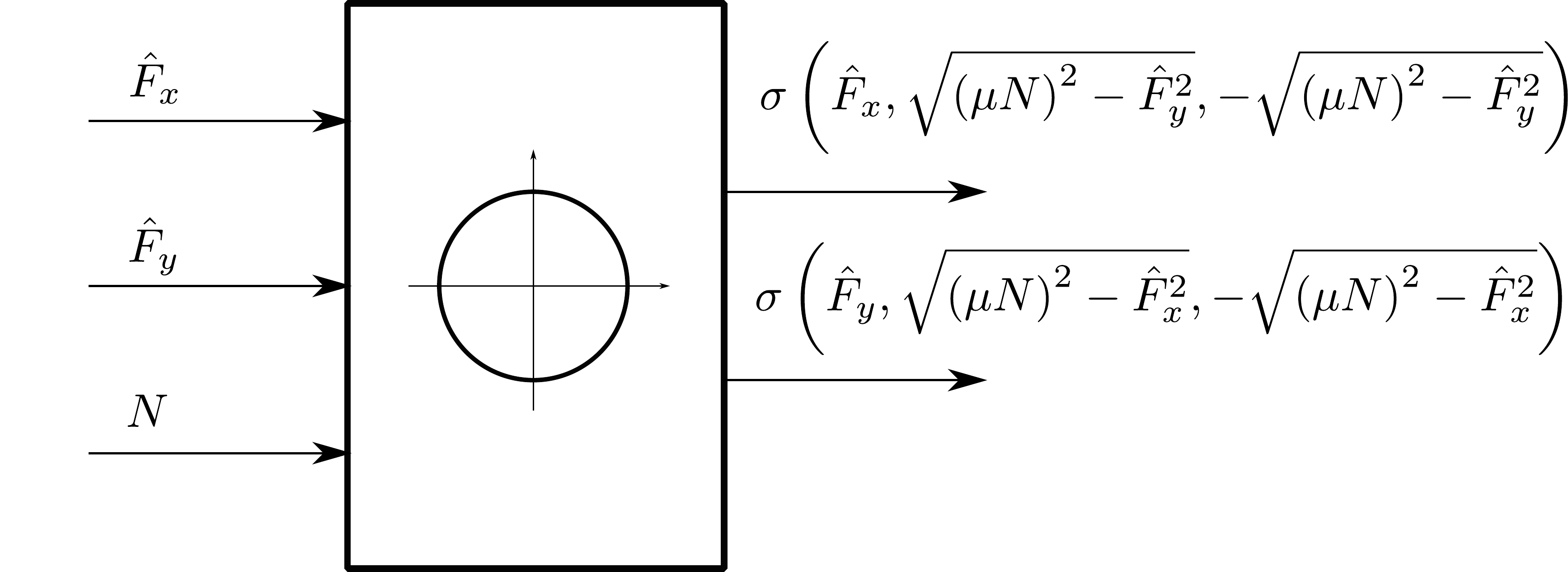}
		\caption{Tire forces saturation.}\label{fig_devMe_saturation_sigma}
	\end{figure}

\begin{figure} [h]
		\centering
		\includegraphics[width = \columnwidth]{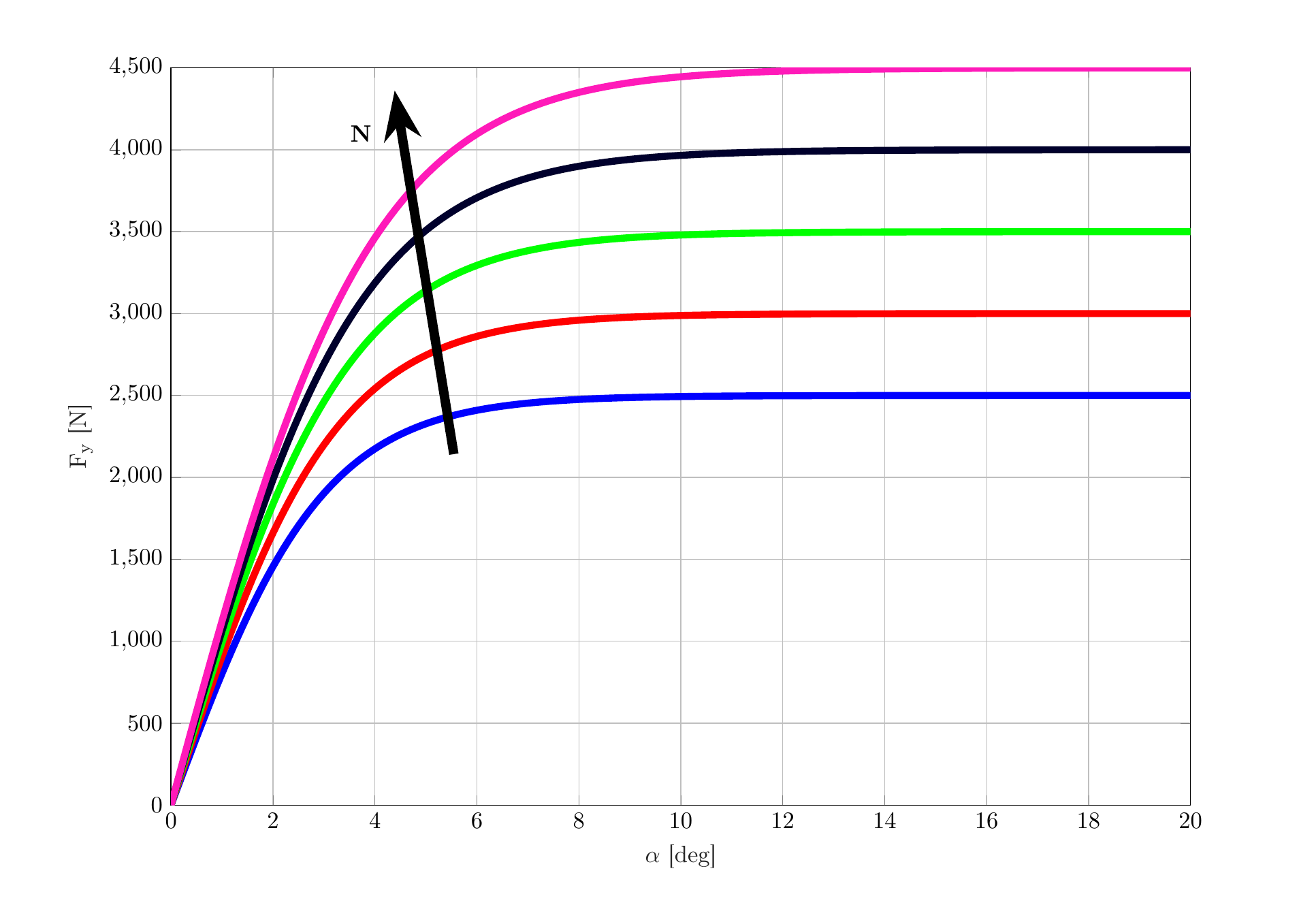}
		\caption{Lateral tire force with respect to $N$.}\label{fig_devMe_tireForces_N_stndAln}
	\end{figure}

	In high dynamics maneuvers, the tires are subject to saturation and combined slip \cite{pacejka_tyre_2012}.
	It is important then to take into consideration these phenomena by using a non-linear model for the tire forces.\\
	The tire model in \cite{kissai_new_2017}, developed for control purposes, is a variation of the \emph{Dugoff's model} \cite{dugoff_tire_1969}.
	It is simpler and requires less parameters than other commonly used models, as the \emph{Pacejka's magic formula} \cite{bakker_tyre_1987} and the \emph{brush model} \cite{pacejka_tyre_2012}, \cite{svendenius_brush_2003}, \cite{villela_nonlinear_2004}.
	We rely on this model since it allows to represent tire forces saturation and combined slip.\\
	Longitudinal and lateral forces are linear respectively on the longitudinal and lateral slip, $\kappa$ and $\alpha$:
	\begin{subequations} \label{eq_devMe_kissaiTireModel}
	\begin{align}
	\hat{F}_x &= c^{*}_{\kappa}\left(\mu, N, \alpha\right) \, \kappa \\
	\hat{F}_y &= c^{*}_{\alpha}\left(\mu, N, \kappa\right) \, \alpha.
	\end{align}
	\end{subequations}
	Longitudinal and lateral slip $\kappa$ and $\alpha$ are defined in \cite{pacejka_tyre_2012}.
	The longitudinal slip, or slip ratio, $\kappa$ is the the difference between wheel longitudinal speed and wheel tangent speed divided by the wheel longitudinal speed.
	The lateral slip angle $\alpha$, shown in Figure \ref{fig_devMe_wheel_slipAngle}, is the difference between the wheel steering angle $\delta$ and the wheel vector speed angle $\theta_w$.
	~\\
	The expressions of $\kappa$ and $\alpha$ for the front and rear wheels are:
	\begin{subequations} \label{eq_devMe_slip_bicycle_3DOF}
	\begin{align}
	\kappa_f &= -\frac{v \, \cos \delta_f + \left(u + \ell_f r\right) \, \sin \delta_f - \omega_{wf} r_{e}}{v \, \cos \delta_f + \left(u + \ell_f r\right) \, \sin \delta_f} \\[8pt]
	\kappa_r &= -\frac{v - \omega_{wr} r_{e}}{v} \\[8pt]
	\alpha_f &= \delta_f - \arctan \left( \frac{u + \ell_f r}{v} \right) \\[8pt]
	\alpha_r &= - \arctan \left( \frac{u - \ell_r r}{v} \right).
	\end{align}
	\end{subequations}
	The coefficients $c^{*}_{\kappa}$ et $c^{*}_{\alpha}$ are:\par 
	\begin{subequations} \label{eq_devMe_kissaiTireModel_nonLinCoeff}
	\scriptsize{
	\begin{align}
	c^{*}_{\kappa}(\mu, N, \alpha) &= \frac{N\, c_{\kappa}\, \mu\, \left(4\, \sqrt{{c_{\alpha}}^2\, {\tan\!\left(\alpha\right)}^2 + {c_{\kappa}}^2\, {\kappa^{*}}^2} + N\, \mu\, \left(\kappa^{*} - 1\right)\right)}{4\, {c_{\alpha}}^2\, {\tan\!\left(\alpha\right)}^2 + 4\, {c_{\kappa}}^2\, {\kappa^{*}}^2}\\[8pt]
	c^{*}_{\alpha}(\mu, N, \kappa) &= \frac{N\, c_{\alpha}\, \mu\, \left(4\, \sqrt{{\alpha^{*}}^2\, {c_{\alpha}}^2 + {c_{\kappa}}^2\, {\kappa}^2} + N\, \mu\, \left(\kappa - 1\right)\right)}{4\, {\alpha^{*}}^2\, {c_{\alpha}}^2 + 4\, {c_{\kappa}}^2\, {\kappa}^2}.
	\end{align}
	}
	\end{subequations}
	$N$ is the normal force on the wheel and $\mu$ the friction coefficient between tire and ground.
	The expressions of $\kappa^{*}$ and $\alpha^{*}$ are:
	\begin{subequations} \label{eq_devMe_kissaiTireModel_coeffStar}
		\begin{align}
		\kappa^{*} &= \frac{N\, \mu\, \left(4\, c_{\kappa} + \sqrt{{N}^2\, {\mu}^2 + 8\, c_{\kappa}\, N\, \mu} + N\, \mu\right)}{8\, {c_{\kappa}}^2} \\
		\alpha^{*} &= \frac{N\, \mu}{2\, c_{\alpha}}.
		\end{align}
	\end{subequations}
	The parameters $c_{\kappa}$ et $c_{\alpha}$ are respectively the longitudinal and lateral tire stiffness.\\
	Here we use a circle to approximate the saturation of the tire forces as in \cite{kritayakirana_autonomous_2012} and  \cite{schuring_model_1996}, thus simplifying the more general ellipse model \cite{pacejka_tyre_2012}, \cite{svendenius_tire_2003}, \cite{wong_theory_2001}.
	The radius of the saturation circle is given by $F_{max} = \mu N$.
	\begin{equation} \label{eq_devMe_kissaiTireModel_saturation}
	\hat{F}_{x}^2 + \hat{F}_{y}^2 \leq \left(\mu N\right)^{2}
	\end{equation}
	~\\
	The logistic function, shown Figure \ref{fig_devMe_logistic_function}, has been used to express the saturation:
	\begin{equation} \label{eq_devMe_logisticFunction}
	\sigma\left( x, u, \ell \right) = \frac{u - \ell}{1+ e^{- \frac{4}{u - \ell} \left(x - \frac{u - \ell}{2}\right)}} + \ell
	\end{equation}
	Where $x$ is the force before saturation, computed as in (\ref{eq_devMe_kissaiTireModel}), $u$ and $\ell$ are the upper and lower saturation bounds.
	A diagram of the tire forces saturation is shown in Figure \ref{fig_devMe_saturation_sigma}.
	~\\
	The normal force $N$ on the tire is present both in the computation of the coefficients $c^{*}_{\kappa}$ et $c^{*}_{\alpha}$, expressed in (\ref{eq_devMe_kissaiTireModel_nonLinCoeff}), and the saturation bound (\ref{eq_devMe_kissaiTireModel_saturation}).
	Considering the importance of $N$ in the computation of the tire forces, as shown in Figure \ref{fig_devMe_tireForces_N_stndAln}, a variable normal force model has then been developed. \\
	We consider the equations of motion of the vehicle neglected movements:
	\begin{subequations} \label{eq_devMe_bicycle_EquationMotion_3DOF_normalForces}
	\begin{align}
	0 &= N_f + N_r - P \label{eq_devMe_bicycle_EquationMotion_3DOF_normalForces_1}\\[8pt]
	0 &= 2 \, \left(F_{yxf} + F_{yyf} + F_{yr}\right) \, h  - 2 \, R_{yxf} \, (h - r_{e}) \label{eq_devMe_bicycle_EquationMotion_3DOF_normalForces_2}\\[8pt]
	\begin{split}
	0 &= -2 \, \left(F_{xxf} - F_{xyf} + F_{xr}\right) \, h \\
	&+ 2 \, \left(R_{xxf} + R_{xr}\right) \, (h - r_{e}) + N_r \, \ell_r - N_f \, \ell_f. \label{eq_devMe_bicycle_EquationMotion_3DOF_normalForces_3}
	\end{split}
	\end{align}
	\end{subequations}
	Equation (\ref{eq_devMe_bicycle_EquationMotion_3DOF_normalForces_1}) describes the vehicle vertical motion.
	Equations (\ref{eq_devMe_bicycle_EquationMotion_3DOF_normalForces_2}) and (\ref{eq_devMe_bicycle_EquationMotion_3DOF_normalForces_3}) describe respectively the roll and pitch motion, i.e. the rotation about the vehicle's $x$ and $y$ axis.
	According to the model definition of Section \ref{subsec_devMe_nonLinearModel}, the vertical, roll and pitch accelerations are equal to zero. \\
	Considering both (\ref{eq_devMe_bicycle_EquationMotion_3DOF_vehicle}) and (\ref{eq_devMe_bicycle_EquationMotion_3DOF_normalForces}), we solve the system of equations with the following unknowns: $N_f$, $N_r$, $F_{xf}$, $F_{xr}$, $F_{yf}$ and $F_{yr}$.
	The solution of $N_f$ and $N_r$ does not depend on the tire model used.
	The expressions obtained for $N_f$ and $N_r$, not shown in this article for the sake of brevity, depend on the vehicle dynamics, in particular the longitudinal acceleration $\dot{v}$.
	~\\
	$N_f$ and $N_r$ are saturated between $0$ and $P_z$ using the logistic function.
	
	\subsection{Model linearization} \label{subsec_devMe_modelLinearization}
\begin{figure} [h]
		\centering
		\includegraphics[width = 0.85\columnwidth]{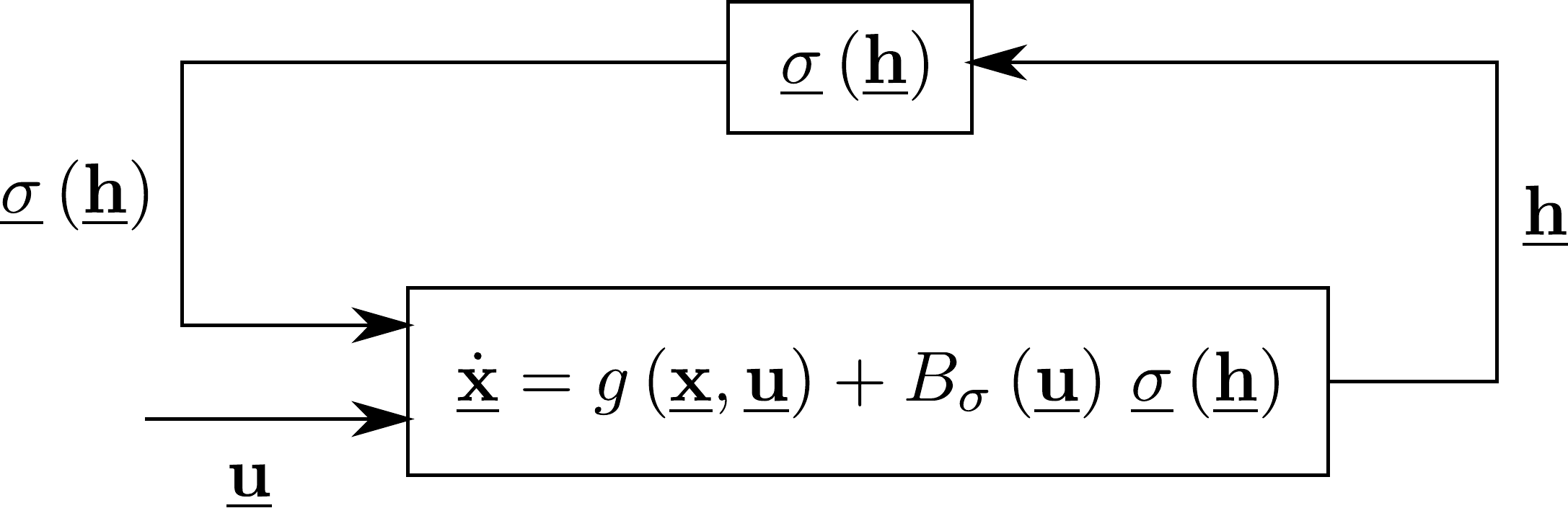}
		\caption{Isolation of the implicit part of (\ref{eq_devMe_bicycle_EquationMotion_3DOF_vehicle}).}\label{fig_devMe_LFT_non_linear}
	\end{figure}

\begin{figure} [h]
		\centering
		\includegraphics[width = \columnwidth]{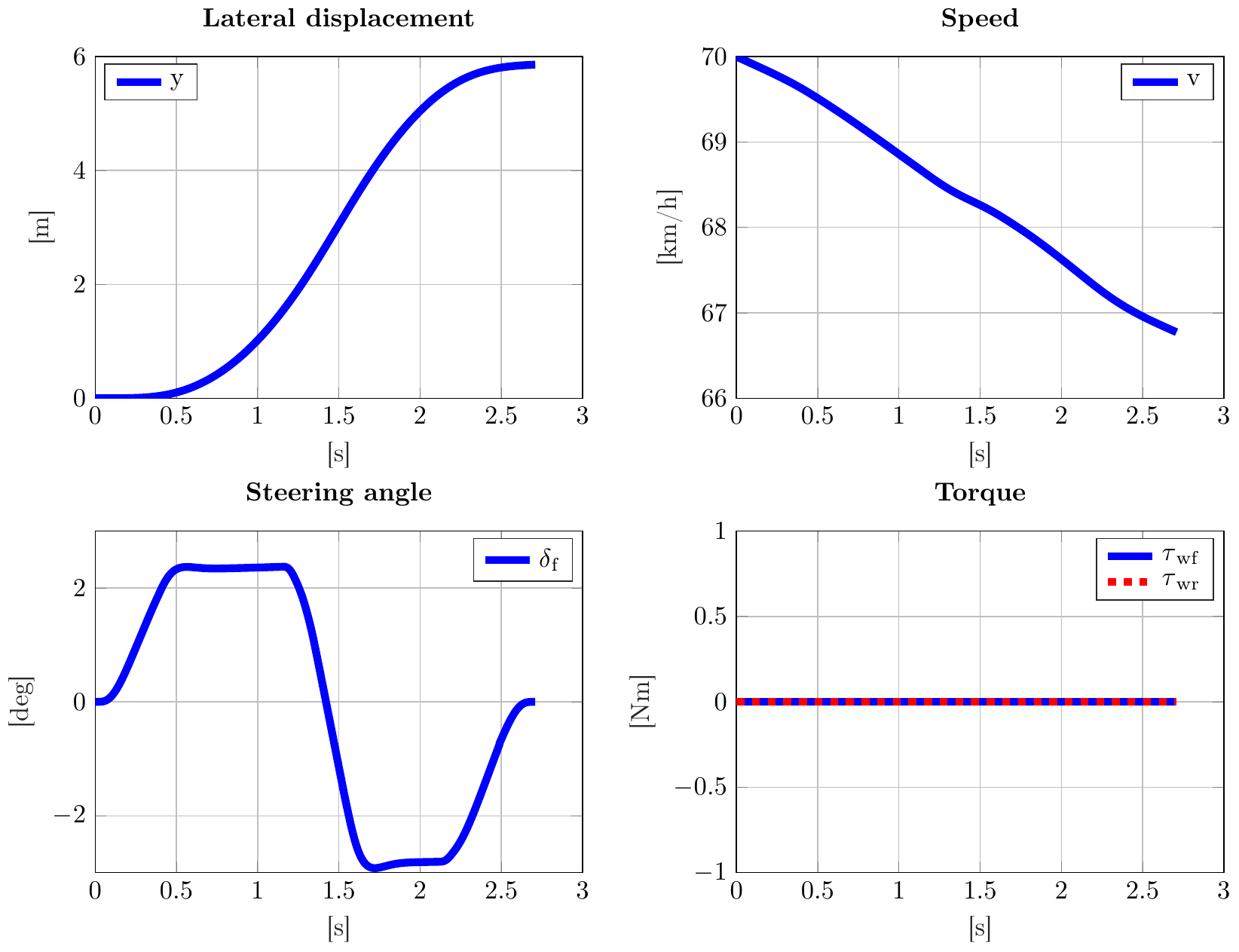}
		\caption{Reference trajectory.}\label{fig_devMe_usCs_info}
	\end{figure}

\begin{figure} [h]
		\centering
		\includegraphics[width = \columnwidth]{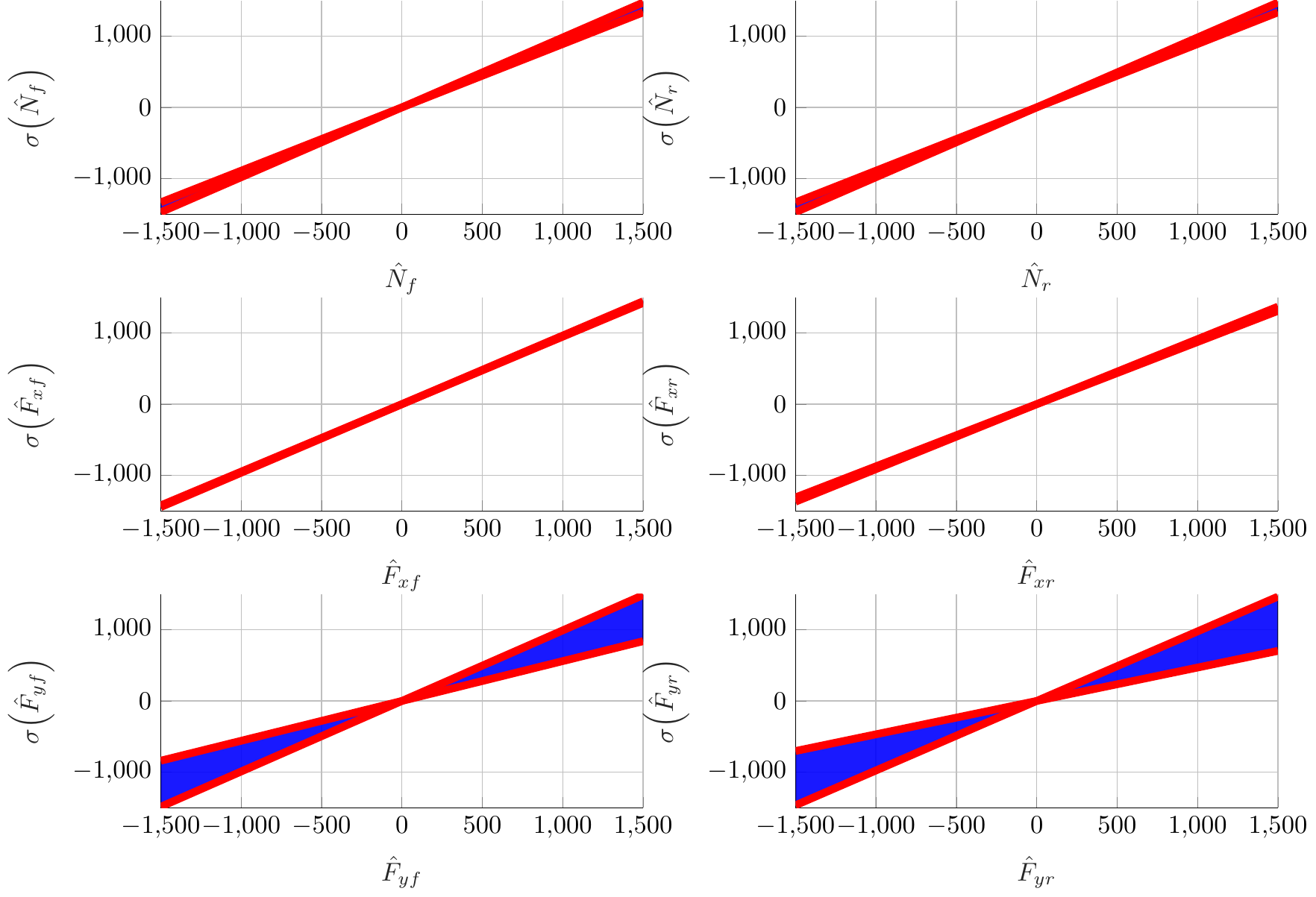}
		\caption{Slope sector for reference trajectory of Figure \ref{fig_devMe_usCs_info}.}\label{fig_devMe_saturation_area_usCs2_70}
	\end{figure}

	The state and input vectors of the model in (\ref{eq_devMe_bicycle_EquationMotion_3DOF_vehicle}) are:
	\begin{equation} \label{eq_devMe_stateInputVectors}
	\begin{aligned}
	\underline{\mathbf{x}} = \left[\begin{array}{c} v\\ u \\ r \\ \omega_{wf}\\ \omega_{wr} \end{array}\right] && && 
	\underline{\mathbf{u}} = \left[\begin{array}{c} \delta_f \\ \tau_{wf} \\ \tau_{wr} \end{array}\right]
	\end{aligned}
	\end{equation}
	The model expressed by (\ref{eq_devMe_bicycle_EquationMotion_3DOF_vehicle}) is in implicit form $f(\dot{\underline{\mathbf{x}}},\underline{\mathbf{x}}, \underline{\mathbf{u}}) = 0$.
	Due to the saturation functions used for the tire forces, it is not possible to express the model in explicit form.\\
	It is instead possible to isolate the implicit part of the model as shown in Figure \ref{fig_devMe_LFT_non_linear}:
	\begin{equation} \label{eq_devMe_BicycleModel_3DOF_implicit_f_LFT}
	\dot{\underline{\mathbf{x}}} = g\left(\underline{\mathbf{x}}, \underline{\mathbf{u}}\right) + B_{\sigma}\left(\underline{\mathbf{u}}\right) \, \underline{\mathbf{\sigma}}\left(\underline{\mathbf{h}}\right)
	\end{equation}
	With:
	\begin{equation} \label{eq_devMe_BicycleModel_3DOF_implicit_sigmaVector}
	\underline{\mathbf{\sigma}}\left(\underline{\mathbf{h}}\right) = \left[\begin{array}{c}
	\sigma \left( h_1\left(\dot{\underline{\mathbf{x}}}, \underline{\mathbf{x}}, \underline{\mathbf{u}}\right)\right) \\
	\sigma \left( h_2\left(\dot{\underline{\mathbf{x}}}, \underline{\mathbf{x}}, \underline{\mathbf{u}}\right)\right) \\
	\sigma \left( h_3\left(\dot{\underline{\mathbf{x}}}, \underline{\mathbf{x}}, \underline{\mathbf{u}}, \underline{\mathbf{\sigma}}\right)\right) \\
	\sigma \left( h_4\left(\dot{\underline{\mathbf{x}}}, \underline{\mathbf{x}}, \underline{\mathbf{u}}, \underline{\mathbf{\sigma}}\right)\right) \\
	\sigma \left( h_5\left(\dot{\underline{\mathbf{x}}}, \underline{\mathbf{x}}, \underline{\mathbf{u}}, \underline{\mathbf{\sigma}}\right)\right) \\
	\sigma \left( h_6\left(\dot{\underline{\mathbf{x}}}, \underline{\mathbf{x}}, \underline{\mathbf{u}}, \underline{\mathbf{\sigma}}\right)\right) \\
	\end{array}\right]
	\end{equation}
	\footnotesize\begin{equation} \label{eq_devMe_BicycleModel_3DOF_implicit_hVector}
	\underline{\mathbf{h}}\left(\dot{\underline{\mathbf{x}}}, \underline{\mathbf{x}}, \underline{\mathbf{u}},  \underline{\mathbf{\sigma}}\right) = \left[\begin{array}{c}
	h_1\left(\dot{\underline{\mathbf{x}}}, \underline{\mathbf{x}}, \underline{\mathbf{u}}\right) \\
	h_2\left(\dot{\underline{\mathbf{x}}}, \underline{\mathbf{x}}, \underline{\mathbf{u}}\right) \\
	h_3\left(\dot{\underline{\mathbf{x}}}, \underline{\mathbf{x}}, \underline{\mathbf{u}}, \underline{\mathbf{\sigma}}\right) \\
	h_4\left(\dot{\underline{\mathbf{x}}}, \underline{\mathbf{x}}, \underline{\mathbf{u}}, \underline{\mathbf{\sigma}}\right) \\
	h_5\left(\dot{\underline{\mathbf{x}}}, \underline{\mathbf{x}}, \underline{\mathbf{u}}, \underline{\mathbf{\sigma}}\right) \\
	h_6\left(\dot{\underline{\mathbf{x}}}, \underline{\mathbf{x}}, \underline{\mathbf{u}}, \underline{\mathbf{\sigma}}\right) \\
	\end{array}\right] = \left[\begin{array}{c}
	\hat{N}_f \\
	\hat{N}_r \\
	\hat{F}_{xf}\left(\sigma\!\left(\hat{N}_f\right)\right) \\
	\hat{F}_{xr}\left(\sigma\!\left(\hat{N}_r\right)\right) \\
	\hat{F}_{yf}\left(\sigma\!\left(\hat{N}_f\right)\right) \\
	\hat{F}_{yr}\left(\sigma\!\left(\hat{N}_r\right)\right)
	\end{array}\right].
	\end{equation} \normalsize
	~\\
	The vector $\underline{\mathbf{h}}\left(\dot{\underline{\mathbf{x}}}, \underline{\mathbf{x}}, \underline{\mathbf{u}},  \underline{\mathbf{\sigma}}\right)$ contains the inputs of the saturation functions.
	The vector $\underline{\mathbf{\sigma}}\left(\underline{\mathbf{h}}\right)$ contains the output of the saturation functions.
	The model non-linearities due to the saturation functions have then been isolated and they result as an input of the model.
	The expressions of $g\left(\underline{\mathbf{x}}, \underline{\mathbf{u}}\right)$ and $B_{\sigma}\left(\underline{\mathbf{u}}\right)$, not shown in this article for the sake of brevity, are derived from (\ref{eq_devMe_bicycle_EquationMotion_3DOF_vehicle}) after extraction of  $\underline{\mathbf{\sigma}}\left(\underline{\mathbf{h}}\right)$. \\
	The state dynamics expressed as in (\ref{eq_devMe_BicycleModel_3DOF_implicit_f_LFT}) can be linearized along a reference trajectory $\left(\dot{\underline{\mathbf{x}}}_{0}(t), \underline{\mathbf{x}}_{0}(t), \underline{\mathbf{u}}_{0}(t), \underline{\mathbf{\sigma}}_{0}(t)\right)$, with $\underline{\mathbf{\sigma}}_{0}(t) = \underline{\mathbf{\sigma}}\!\left(\underline{\mathbf{h}}\!\left(\dot{\underline{\mathbf{x}}}_{0}(t), \underline{\mathbf{x}}_{0}(t), \underline{\mathbf{u}}_{0}(t) \right) \right)$.
	For simplicity, in the following the reference trajectory will be indicated with $\mathscr{T}(t)$. The vector $\underline{\mathbf{\sigma}}\left(\underline{\mathbf{h}}\right)$ is also non linear with respect to $\underline{\mathbf{h}}\left(\dot{\underline{\mathbf{x}}}, \underline{\mathbf{x}}, \underline{\mathbf{u}}, \underline{\mathbf{\sigma}}\right)$.
	It is linearized too along the trajectory $\mathscr{T}(t)$. \par
	\begin{subequations}\label{eq_devMe_BicycleModel_3DOF_LFT_linearisation_all}
	\small
	\begin{gather}
	\Delta \dot{\underline{\mathbf{x}}} = A(t) \Delta \underline{\mathbf{x}} + B(t) \Delta \underline{\mathbf{u}} + B_{\sigma}(t) \Delta \underline{\mathbf{\sigma}} \label{eq_devMe_BicycleModel_3DOF_LFT_linearisation_Delta_x}\\
	\Delta \underline{\mathbf{h}} = C(t) \Delta \underline{\mathbf{x}} + D(t) \Delta \underline{\mathbf{u}} + D_{\sigma}(t) \Delta \underline{\mathbf{\sigma}} \label{eq_devMe_BicycleModel_3DOF_LFT_linearisation_Delta_h} \\[8pt]
	A(t) = \left.\frac{\partial g\left(\underline{\mathbf{x}}, \underline{\mathbf{u}}\right)}{\partial \underline{\mathbf{x}}}\right|_{\mathscr{T}(t)} \\[8pt]
	B(t) = \left.\frac{\partial g\left(\underline{\mathbf{x}}, \underline{\mathbf{u}}\right)}{\partial \underline{\mathbf{u}}}\right|_{\mathscr{T}(t)} + \left.\frac{\partial \, \left(B_{\sigma}\left(\underline{\mathbf{u}}\right) \, \underline{\mathbf{\sigma}}\left(\underline{\mathbf{h}}\right)\right)}{\partial \underline{\mathbf{u}}}\right|_{\mathscr{T}(t)} \\[8pt]
	B_{\sigma}(t) = \left.B_{\sigma}\left(\underline{\mathbf{u}}\right)\right|_{\mathscr{T}(t)} \\[8pt]
	C(t) = \left.\frac{\partial \underline{\mathbf{h}}\left(\dot{\underline{\mathbf{x}}}, \underline{\mathbf{x}}, \underline{\mathbf{u}}, \underline{\mathbf{\sigma}}\right)}{\partial \dot{\underline{\mathbf{x}}}}\right|_{\mathscr{T}(t)} \, A(t) + \left.\frac{\partial \underline{\mathbf{h}}\left(\dot{\underline{\mathbf{x}}}, \underline{\mathbf{x}}, \underline{\mathbf{u}}, \underline{\mathbf{\sigma}}\right)}{\partial \underline{\mathbf{x}}}\right|_{\mathscr{T}(t)}
	\\[8pt]
	D(t) = \left.\frac{\partial \underline{\mathbf{h}}\left(\dot{\underline{\mathbf{x}}}, \underline{\mathbf{x}}, \underline{\mathbf{u}}, \underline{\mathbf{\sigma}}\right)}{\partial \dot{\underline{\mathbf{x}}}}\right|_{\mathscr{T}(t)} \, B(t) + \left.\frac{\partial \underline{\mathbf{h}}\left(\dot{\underline{\mathbf{x}}}, \underline{\mathbf{x}}, \underline{\mathbf{u}}, \underline{\mathbf{\sigma}}\right)}{\partial \underline{\mathbf{u}}}\right|_{\mathscr{T}(t)}
	\\[8pt]
	D_{\sigma}(t) = \left.\frac{\partial \underline{\mathbf{h}}\left(\dot{\underline{\mathbf{x}}}, \underline{\mathbf{x}}, \underline{\mathbf{u}}, \underline{\mathbf{\sigma}}\right)}{\partial \dot{\underline{\mathbf{x}}}}\right|_{\mathscr{T}(t)} \, B_{\sigma}(t) + \left.\frac{\partial \underline{\mathbf{h}}\left(\dot{\underline{\mathbf{x}}}, \underline{\mathbf{x}}, \underline{\mathbf{u}}, \underline{\mathbf{\sigma}}\right)}{\partial \underline{\mathbf{\sigma}}}\right|_{\mathscr{T}(t)}
	\end{gather}
	\normalsize
	\end{subequations}
	In this work the reference trajectory of Figure \ref{fig_devMe_usCs_info} has been considered.
	The trajectory has a steering to the left followed by an opposite steering to the right.
	That results in a vehicle's lateral displacement of approximately 6m.
	No torque is applied to the wheels.
	This trajectory corresponds to a two lane change maneuver at 70km/h. \\ 
	The saturation functions can be approximated as sector-bounded non-linearities.
	This means that they are approximated to a line, whose slope varies according to the system working point.
	Figure \ref{fig_devMe_saturation_area_usCs2_70} shows the slope sector for each element of $\underline{\mathbf{\sigma}}\left(\underline{\mathbf{h}}\right)$ along the reference trajectory.
	Since in the reference trajectory the wheel torques are equal to zero, the longitudinal slip ratios are small, hence the longitudinal tire forces are small too.
	As a consequence they are far from their saturation bounds and their slope sector is small.
	The absence of torque on the wheel causes a small longitudinal acceleration $\dot{v}$.
	Since the wheel normal forces vary with $\dot{v}$, as discussed in Section \ref{subsec_devMe_nonLinearModel}, the load transfer is small. 
	Therefore their slope sector is small too.
	On the contrary the slope sector of the lateral forces is large.
	Indeed the input steering angle $\delta_f$ produces large lateral forces, which move closer to their saturation bounds.
	~\\
	It is then possible to consider $\underline{\mathbf{\sigma}}\left(\underline{\mathbf{h}}\right) \approx K_{\sigma} \underline{\mathbf{h}}$. $K_{\sigma}$ is a diagonal matrix, with the slope of the saturation functions along the diagonal.
	Applying this approximation in the expression of $\Delta \underline{\mathbf{h}}$ in (\ref{eq_devMe_BicycleModel_3DOF_LFT_linearisation_Delta_h}) we obtain:
	\begin{equation*}
	\begin{aligned}
	\Delta \underline{\mathbf{h}} =& C(t) \Delta \underline{\mathbf{x}} + D(t) \Delta \underline{\mathbf{u}} + D_{\sigma}(t) \, K_{\sigma} \Delta \underline{\mathbf{h}} \\[8pt]
	\Delta \underline{\mathbf{h}} =& \left(I - D_{\sigma}(t) \, K_{\sigma} \right)^{-1} \, C(t) \Delta \underline{\mathbf{x}} \\
	&+ \left(I - D_{\sigma}(t) \, K_{\sigma} \right)^{-1} \, D(t) \Delta \underline{\mathbf{u}}.
	\end{aligned}
	\end{equation*}
	Replacing this expression of $\Delta \underline{\mathbf{h}}$ in (\ref{eq_devMe_BicycleModel_3DOF_LFT_linearisation_Delta_x}), the linearized state dynamics with the sector-bounded approximation of the saturation functions is:
	\begin{subequations}\label{eq_devMe_BicycleModel_3DOF_LFT_linearisation_der_x_satApprox1}
	\begin{gather}
	\Delta \dot{\underline{\mathbf{x}}} = \widetilde{A}(t) \Delta \underline{\mathbf{x}} + \widetilde{B}(t) \Delta \underline{\mathbf{u}} \\[8pt]
	\widetilde{A}(t) = A(t) + B_{\sigma}(t) \, K_{\sigma} \, \left(I - D_{\sigma}(t) \, K_{\sigma} \right)^{-1} \, C(t) \label{eq_devMe_BicycleModel_3DOF_LFT_linearisation_der_x_satApprox1_A} \\
	\widetilde{B}(t) = B(t) + B_{\sigma}(t) \, K_{\sigma} \, \left(I - D_{\sigma}(t) \, K_{\sigma} \right)^{-1} \, D(t). \label{eq_devMe_BicycleModel_3DOF_LFT_linearisation_der_x_satApprox1_B}
	\end{gather}
	\end{subequations}

	\subsection{LPV polytopic model for control synthesis} \label{subsec_devMe_LPVpolytopic}
	The matrices $\widetilde{A}(t)$ and $\widetilde{B}(t)$ in (\ref{eq_devMe_BicycleModel_3DOF_LFT_linearisation_der_x_satApprox1_A}) and (\ref{eq_devMe_BicycleModel_3DOF_LFT_linearisation_der_x_satApprox1_B}) are time-varying.\\
	In the following, the elements of the matrix $K_{\sigma}$ will be considered as the average slope coefficients of the slope sector shown in Figure \ref{fig_devMe_saturation_area_usCs2_70}. \\
	It is possible to consider the elements of $\widetilde{A}(t)$ and $\widetilde{B}(t)$ that vary with time as varying parameters.
	We call $\theta_i, \: i = 1, \, ... \, , q$ these parameters.
	$\underline{\mathbf{\Theta}}$ is the vector containing all the parameters.
	It is then possible to express the linearized model in (\ref{eq_devMe_BicycleModel_3DOF_LFT_linearisation_der_x_satApprox1}) as an LPV model, with the state matrix $A$ and the input matrix $B$ affinely dependent on the varying parameters. \\
	Along the reference trajectory shown in Figure \ref{fig_devMe_usCs_info}, a total of 23 elements of $\widetilde{A}(t)$ and $\widetilde{B}(t)$ vary with time.
	However several of these elements have a small variation.
	It is possible to consider just the 6 elements of $\widetilde{A}(t)$ and $\widetilde{B}(t)$ that vary the most to constitute the $\underline{\mathbf{\Theta}}$ vector.
	The range of each parameter $\theta_i, \: i = 1, \, ... \, , 6$ is between a minimum $\underline{\theta}_i$ and a maximum $\overline{\theta}_i$.
	$\underline{\theta}_i$ and $\overline{\theta}_i$ are the minimum and maximum values of $\theta_i$ along the reference trajectory $\mathscr{T}(t)$.\\
	For the control synthesis we add the longitudinal, lateral and yaw angle errors with respect to $\mathscr{T}(t)$ to the state of the model. They are defined as:
	\begin{subequations}\label{eq_devMe_error}
		\begin{align}
		x_L &= x - x_0 \\
		y_L &= y - y_0 \\
		\Delta \psi &= \psi - \psi_0.
		\end{align}
	\end{subequations}
	$x_0$, $y_0$ and $\psi_0$ are the vehicle's absolute position and attitude along the reference trajectory $\mathscr{T}(t)$. \\
	The vehicle's position dynamics in (\ref{eq_devMe_bicycle_posDynamics_1}) and (\ref{eq_devMe_bicycle_posDynamics_2}) can be simplified neglecting the vehicle's lateral speed $u$:
	\begin{subequations}\label{eq_devMe_kinematicModel}
		\begin{align}
		\dot{x} &= v \, \cos \psi \\
		\dot{y} &= v \, \sin \psi.
		\end{align}
	\end{subequations}
	Considering that $v = v_0 + \Delta v$, $u = u_0 + \Delta u$ and $\psi = \psi_0 + \Delta \psi$ and that the $\Delta \psi$ is small along the reference trajectory, the dynamics of the error variables are:
	\begin{subequations}\label{eq_devMe_errorDynamics}
		\begin{align}
		\dot{x}_L &= \Delta v \\
		\dot{y}_L &= \Delta v \, \psi_0 + \Delta \psi \, v_0 \, \cos \psi_0 \\
		\Delta \dot{\psi} &= \Delta r.
		\end{align}
	\end{subequations}
	The expressions $\psi_0$ and $v_0 \, \cos \psi_0$ are added to $\underline{\mathbf{\Theta}}$. \\
	It is finally possible to define the polytope as the hypercube of dimension 8, each dimension being the parameter $\theta_i, \: i = 1, \, ... \, , 8$, between $\underline{\theta}_i$ and $\overline{\theta}_i$.
	The polytopic model is defined at its $N = 2^8$ vertices:
	\begin{equation}\label{eq_devMe_LPVpolytopic_ss}
	S_{i} = \left( \begin{array}{cc} A_{i} & B_{i} \\ C_{i} & D_{i} \end{array} \right) \quad i = 1, \, ... \,, N.
	\end{equation} 
	
	\section{Main results} \label{sec_devMe_main}
	A static state feedback controller has been designed starting from the polytopic LPV model described in Section \ref{subsec_devMe_LPVpolytopic}.
	A stabilizing state feedback for (\ref{eq_devMe_LPVpolytopic_ss}) can be computed solving a problem under LMI constraints, as shown in \cite{blanchini_set-theoretic_2015}.
	The system (\ref{eq_devMe_LPVpolytopic_ss}) is stabilizable by a feedback controller $K = R \, Q^{-1}$ with a specific level of contractivity $\beta > 0$ if there exist a matrix $Q = Q^{T} > 0$ and a matrix $R$ such that: \par
	\small\begin{equation}\label{eq_devMe_staticFeedback_Blanchini}
	Q A_{i}^{T}+A_{i} Q+R^{T} B_{i}^{T}+B_{i} R + 2 \beta Q < 0, \quad i=1, \ldots, N.
	\end{equation}\normalsize
	The LMI problem expressed in (\ref{eq_devMe_staticFeedback_Blanchini}) allows to specify, as performance criteria, just the level of contractivity. \\
	In order to specify other constraints on the poles of the closed loop, it is possible to define an LMI region as in \cite{chilali_robust_1999}. \\
	An LMI region $\mathcal{D}$ of order $s$ is a subset of the complex plane defined as:
	\begin{equation}\label{eq_devMe_LMI_region}
	\mathcal{D}=\left\{z \in \mathbb{C}: \quad L+z M+\bar{z} M^{T}<0\right\}
	\end{equation}
	With $L = L^{T} \in \mathbb{R}^{s \times s}$ and $M \in \mathbb{R}^{s \times s}$. \\
	The system (\ref{eq_devMe_LPVpolytopic_ss}) is $\mathcal{D}$\textit{-stabilizable} by a state feedback controller $K = R \, Q^{-1}$ if there exist a matrix $Q = Q^{T} > 0$ and a matrix $R$ such that: \par
\begin{figure} [h]
		\centering
		\includegraphics[width = 0.55\columnwidth]{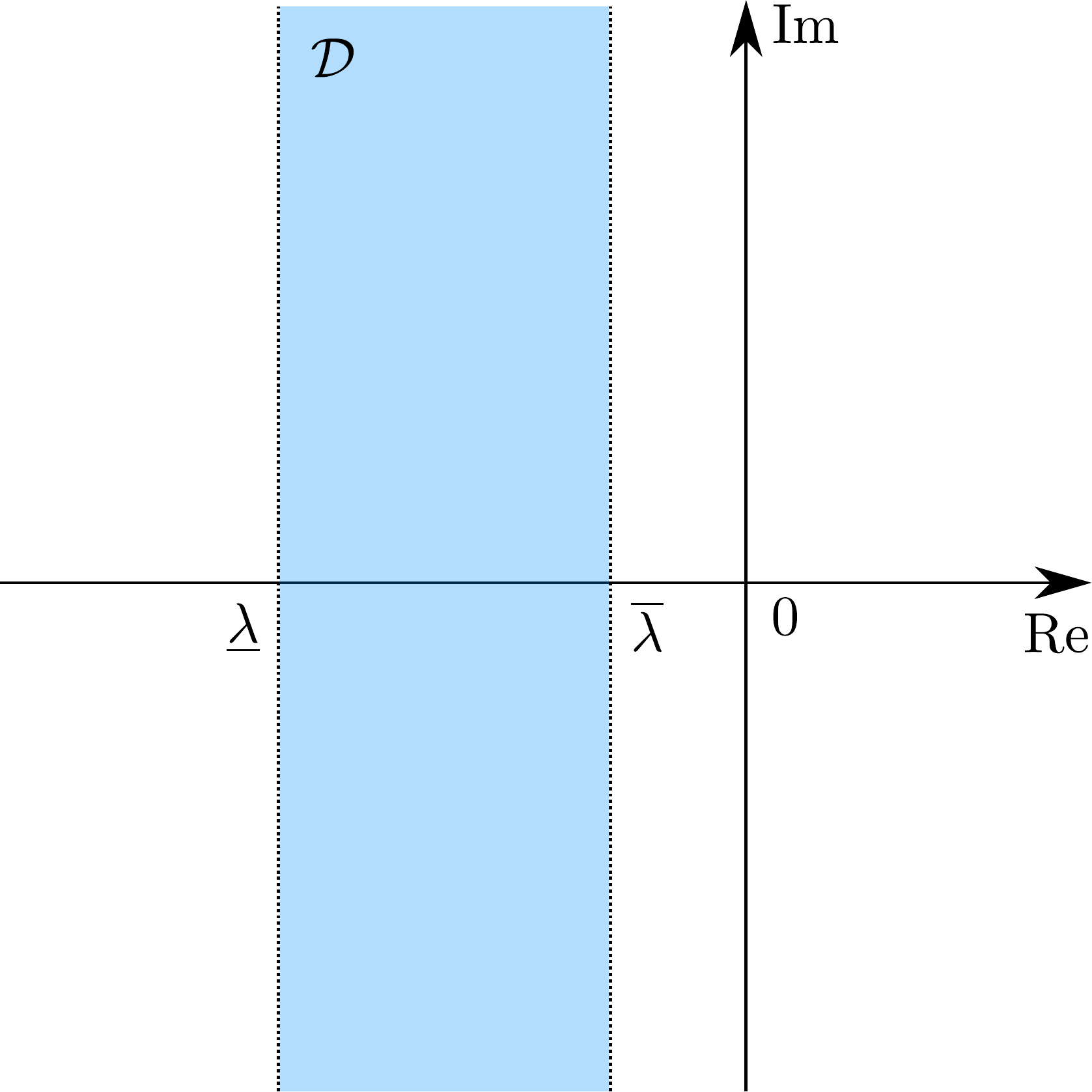}
		\caption{LMI region for control synthesis.}\label{fig_devMe_LMI_region}
	\end{figure}

	\footnotesize\begin{equation}\label{eq_devMe_Dstabilizable}
	\begin{split}
	L& \otimes Q+M \otimes\left(A_i Q\right)+M^{T} \otimes\left(Q A_i^{T}\right)\\
	&+ M \otimes\left(B_i L\right)+M^{T} \otimes\left(L^{T} B_i^{T}\right)<0
	\end{split} \quad \quad i=1, \ldots, N.
	\end{equation}\normalsize
	The controller presented in this paper has been developed using this last method. \\
	The LMI region $\mathcal{D}$ has been defined in order to restrict the dynamics of the closed loop between $\overline{\lambda} = -2$ and $\underline{\lambda} = -40$, as shown in Figure \ref{fig_devMe_LMI_region}.
	This choice guarantees a certain level of contractivity and limits fast dynamics.
	
	\section{Numerical results} \label{sec_devMe_numericalResults} 
\begin{figure} [h!]
		\centering
		\includegraphics[width = 1\columnwidth]{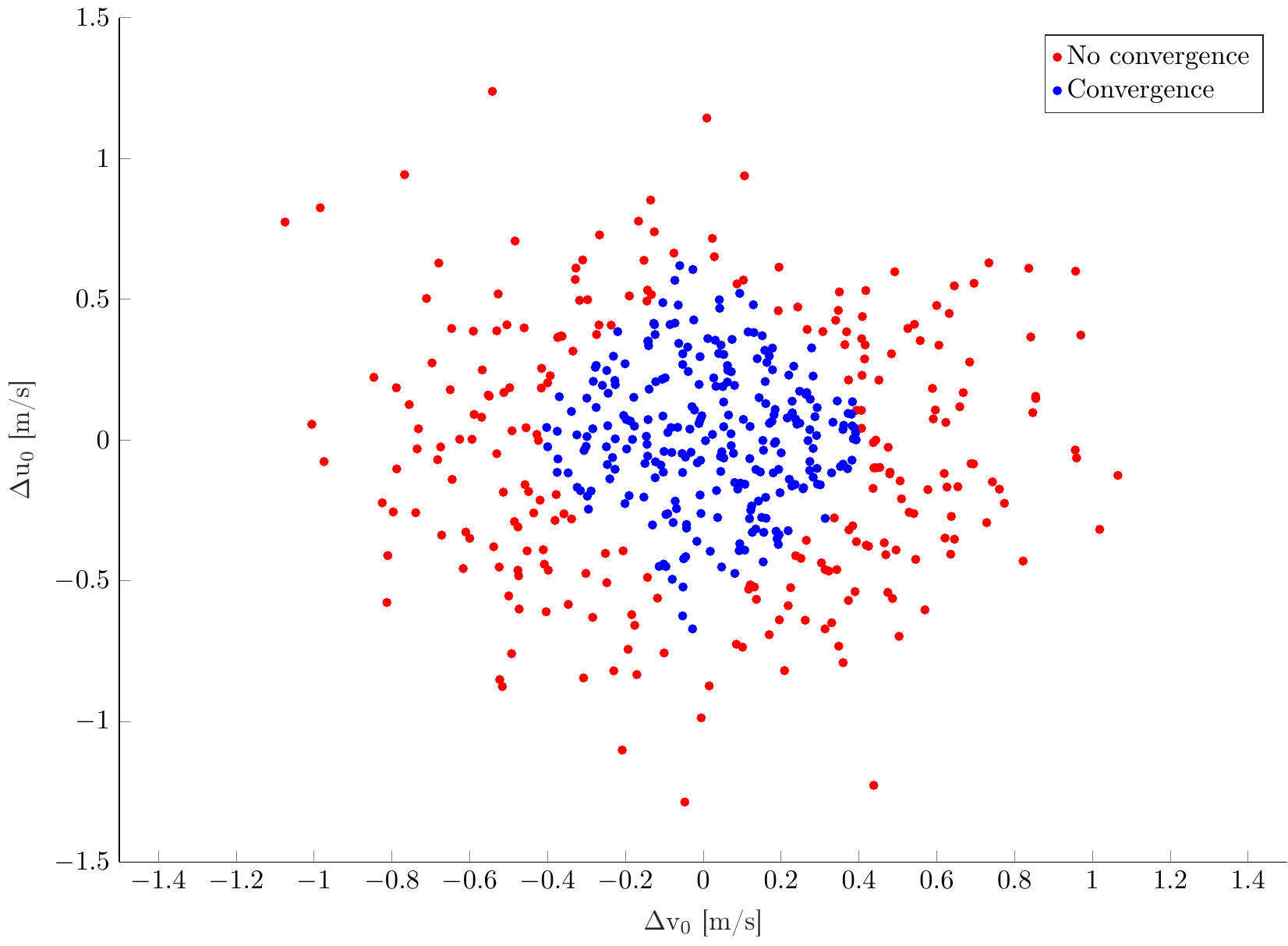}
		\caption{Empirical region of attraction.}\label{fig_devMe_v_u_scatter}
	\end{figure}

\begin{figure} [h!]
		\centering
		\includegraphics[width = 1\columnwidth]{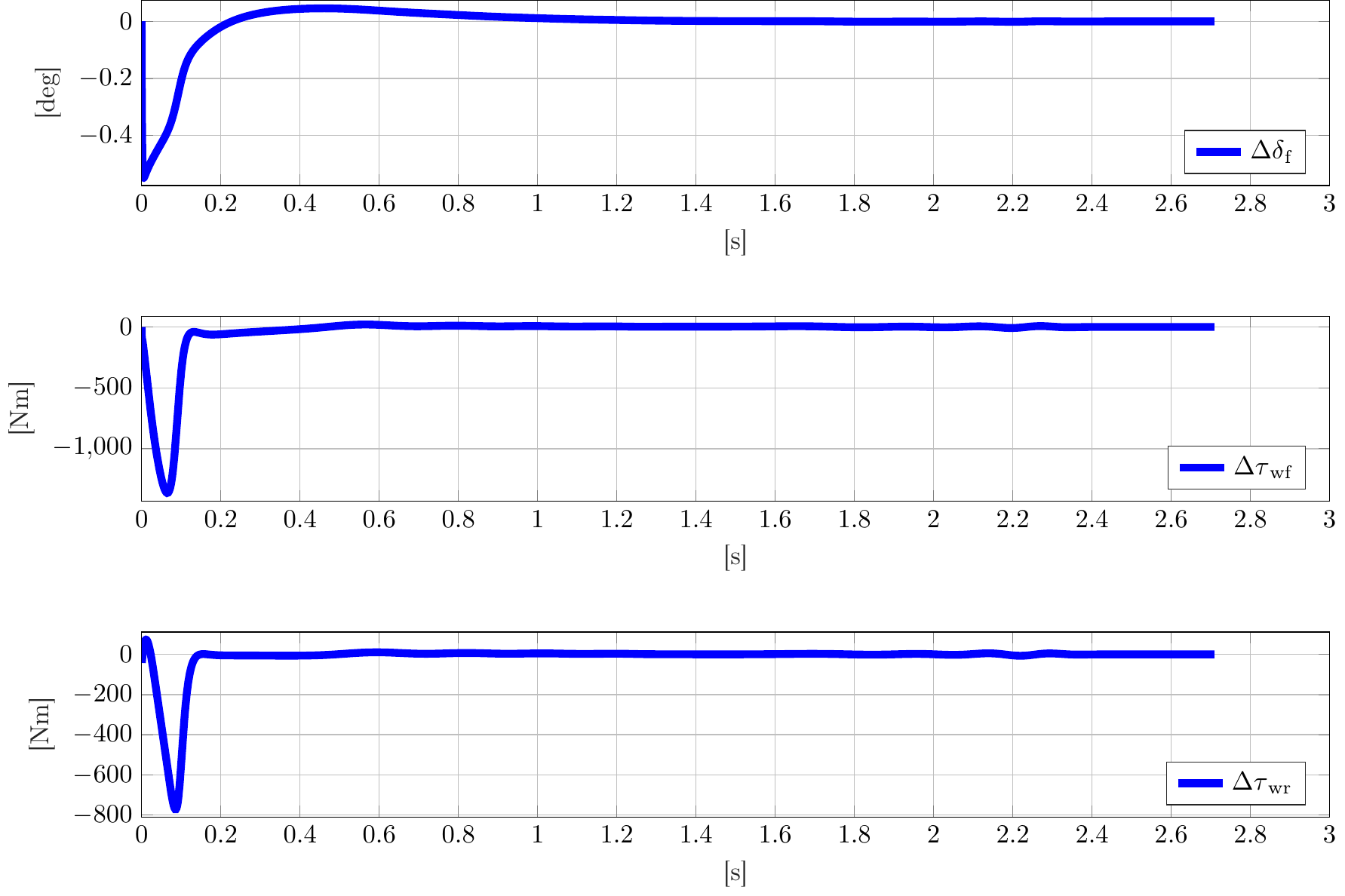}
		\caption{Commands computed by the controller.}\label{fig_devMe_cmnd_v03_u03}
	\end{figure}

\begin{figure} [h!]
		\centering
		\includegraphics[width = 1\columnwidth]{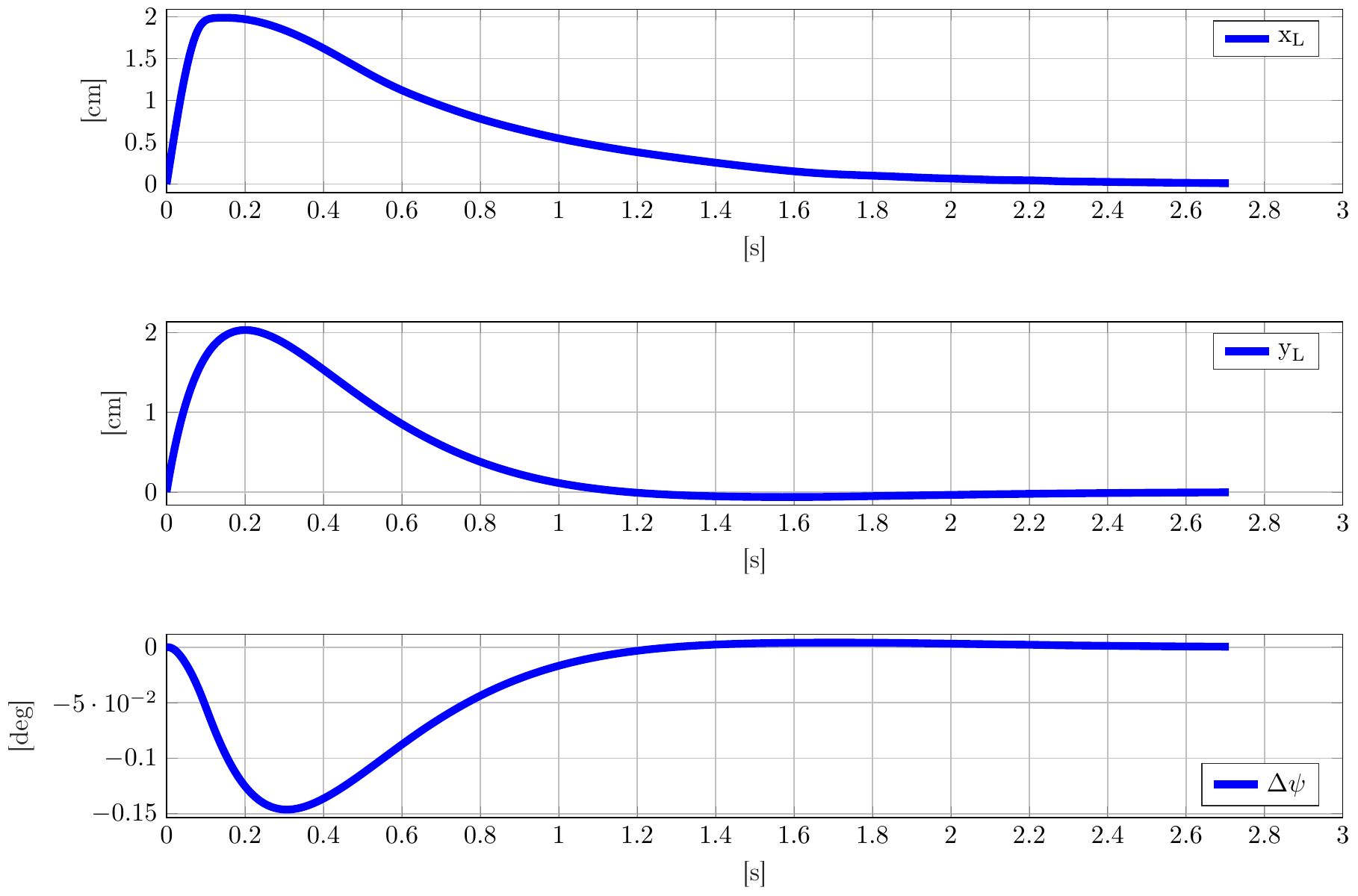}
		\caption{Vehicle's position and attitude error.}\label{fig_devMe_error_v03_u03}
	\end{figure}

\begin{figure} [h!]
		\centering
		\includegraphics[width = 1\columnwidth]{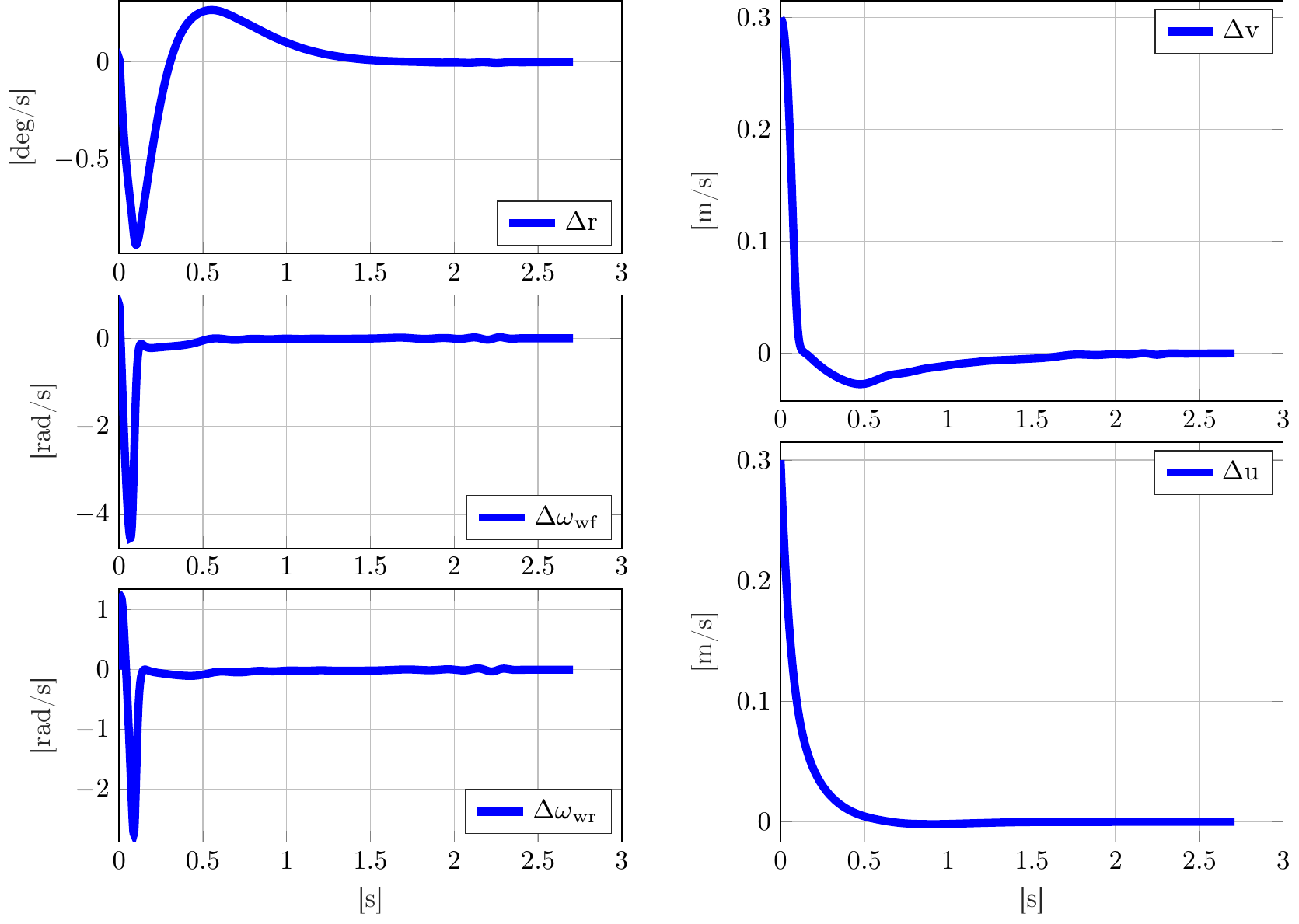}
		\caption{Vehicle's dynamics error.}\label{fig_devMe_state_v03_u03}
	\end{figure}

	In this Section some numerical results of the controller are shown. \\
	The controller has been tested on Simulink on the non-linear bicycle model described in Section \ref{subsec_devMe_nonLinearModel}.
	The trajectory considered is the reference trajectory shown in Section \ref{subsec_devMe_modelLinearization}, Figure \ref{fig_devMe_usCs_info}. \\
	Several simulations have been run, with a set of different initial conditions for $\Delta v$ and $\Delta u$. 
	The purpose is to establish empirically the region of initial conditions for which the controller is able to stabilize the system and converge to the reference trajectory. 
	Figure \ref{fig_devMe_v_u_scatter} shows the empirical region of attraction of the controller for a non-zero initial conditions for the variables $\Delta v$ and $\Delta u$.\\
	The capacity of the controller to stabilize the system and bring back the vehicle along the reference trajectory is limited for $\mid \Delta v \left( 0 \right) \mid \lesssim 0.4$m/s and $\mid \Delta u \left( 0 \right) \mid \lesssim 0.64$m/s. \\
	A case with initial conditions of $\Delta v \left( 0 \right) = 0.3$m/s and $\Delta u \left( 0 \right) = 0.3$m/s is shown more in detail. \\
	Figure \ref{fig_devMe_cmnd_v03_u03} shows the commands computed by the controller.
	Figure \ref{fig_devMe_error_v03_u03} shows the vehicle's position and attitude error with respect to the reference trajectory.
	Finally Figure \ref{fig_devMe_state_v03_u03} shows the vehicle's dynamics error with respect to the reference trajectory.
	Both the dynamics and position errors converge to zero, i.e. the vehicle is brought back to the reference trajectory.
	~\\
	It has been observed that in the cases in which the closed loop diverges, the torques computed by the controller are too large, outside the limits of what could be possibly realized on a vehicle.
	
	\section{Conclusion and perspectives} \label{sec_devMe_conclusion}
	The results of Section \ref{sec_devMe_numericalResults} show that the LPV model developed in Section \ref{sec_devMe_preliminaries} can be successfully used to develop a controller that stabilizes the vehicle along a collision avoidance maneuver.
	Even if the set of initial conditions for which the controller is able to stabilize the system is limited, the LPV systems framework offers many possibilities to improve the controller.
	The LPV model could be used for the synthesis of a dynamic feedback controller, as the one proposed in \cite{apkarian_self-scheduled_1995}. \\
	The tire forces will also be studied in the absolute stability context, that is, considering the whole sector for the nonlinearity description instead of considering a single slope approximation.
	
	\bibliographystyle{IEEEtran}
	\bibliography{IEEEabrv,mylib_20210126}
\end{document}